\documentclass{amsart}
\usepackage{amsfonts,amssymb,amscd,amsmath,enumerate,verbatim}
\usepackage{xypic}
\xyoption{curve}

\usepackage{hyperref}

\theoremstyle{plain}

\newtheorem{theorem}{Theorem}[section]
\newtheorem{proposition}[theorem]{Proposition}
\newtheorem{lemma}[theorem]{Lemma}
\newtheorem{corollary}[theorem]{Corollary}

\theoremstyle{definition}

\newtheorem{examples}[theorem]{Examples}

\theoremstyle{remark}
\newtheorem{remark}[theorem]{Remark}

\newtheorem{chunk}[theorem]{}

\numberwithin{equation}{theorem}

\newcommand{\lra}{\longrightarrow}

\newcommand{\xra}{\xrightarrow}

\newcommand{\up}[1]{{{}^{#1}\!}}

\newcommand{\wh}{\widehat}

\newcommand{\les}{{\scriptscriptstyle\leqslant}}

\newcommand{\dcat}[1]{{\mathsf D}(#1)}
\newcommand{\dcatb}[1]{{\mathsf{D}^{\mathsf f}_{\mathsf b}}(#1)}
\newcommand{\dcatn}[1]{{\mathsf{D}^{\mathsf f}_{\!{\scriptscriptstyle\mathsf -}}}(#1)}
\newcommand{\dcatp}[1]{{\mathsf{D}^{\mathsf f}_{\!{\scriptscriptstyle\mathsf +}}}(#1)}
\newcommand{\dtensor}[1]{\otimes_{#1}^{\mathsf {L}}}
\newcommand{\koszul}[2]{{{K}[#1;#2]}}
\newcommand{\Rhom}[3]{{\mathsf{R}\!\Hom}_{#1}(#2,#3)}

\newcommand{\shift}{{\sf\Sigma}}

\newcommand\dd{\partial}
\newcommand{\hh}[1]{\operatorname{H}(#1)}
\newcommand\HH{\operatorname{H}}

\newcommand\ZZ{\operatorname{Z}}
\newcommand\Tor{\operatorname{Tor}}
\newcommand\Ext{\operatorname{Ext}}
\newcommand\Hom{\operatorname{Hom}}

\newcommand\Spec{\operatorname{Spec}}
\newcommand\Supp{\operatorname{Supp}}

\newcommand\codim{\operatorname{codim}}

\newcommand\edim{\operatorname{edim}}

\newcommand{\lol}{{\ell\ell}}
\newcommand{\dll}[2]{{\ell\ell_{\mathsf{D}(#1)}{#2}}}

\newcommand\rank{\operatorname{rank}}

\newcommand{\Ann}{\operatorname{Ann}}

\newcommand\idmap{\operatorname{id}}

\newcommand{\rP}[2]{P^{#1}_{#2}(t)}
\newcommand{\rI}[2]{I_{#1}^{#2}(t)}

\newcommand{\betti}[3]{\beta_{#1}^{#2}(#3)}
\newcommand{\bass}[3]{\mu^{#1}_{#2}(#3)}

\newcommand\cxy{\operatorname{cx}}
\newcommand\curv{\operatorname{curv}}

\newcommand\injcxy{\operatorname{inj\,cx}}
\newcommand\injcurv{\operatorname{inj\,curv}}

\newcommand\fm{{\mathfrak m}}
\newcommand\fn{{\mathfrak n}}

\newcommand\bsx{{\boldsymbol x}}
\newcommand\bsy{{\boldsymbol y}}
\newcommand\bsz{{\boldsymbol z}}

\newcommand\BN{{\mathbb N}}
\newcommand\BQ{{\mathbb Q}}
\newcommand\BR{{\mathbb R}}
\newcommand\BZ{{\mathbb Z}}

\newcommand{\vf}{{\varphi}}

\begin{document}

\title[Contracting endomorphisms]{Homological invariants of modules\\ over contracting endomorphisms}

\date{20th May 2011}

\author[L.~L.~Avramov]{Luchezar L.~Avramov}
\address{Department of Mathematics,
University of Nebraska, Lincoln, NE 68588, U.S.A.}
\email{avramov@math.unl.edu}

\author[M.~Hochster]{Melvin Hochster}
\address{Department of Mathematics,
University of Michigan, Ann Arbor, MI 48109, U.S.A.}
\email{hochster@umich.edu}

\author[S.~B.~Iyengar]{Srikanth B.~Iyengar}
\address{Department of Mathematics,
University of Nebraska, Lincoln, NE 68588, U.S.A.}
\email{iyengar@math.unl.edu}

\author[Y.~Yao]{Yongwei Yao}
\address{Department of Mathematics and Statistics,
Georgia State University, Atlanta, GA 30303, U.S.A.}
\email{yyao@gsu.edu}

\thanks
{Partly supported by NSF grants DMS-0803082 (LLA), DMS-0400633 and DMS-0901145 (MH), 
DMS-0903493 (SBI), and  DMS-0700554 (YY)}

\begin{abstract}
It is proved that when $R$ is a local ring of positive characteristic, $\phi\colon R\to R$
is its Frobenius endomorphism, and some non-zero finite $R$-module has finite flat 
dimension or finite injective dimension for the $R$-module structure induced through $\phi$, 
then $R$ is regular.  This broad generalization of Kunz's characterization of regularity in 
positive characteristic is deduced from a theorem concerning a local ring $R$ with 
residue field of $k$ of arbitrary characteristic:  If $\phi$ is a contracting endomorphism of  $R$,
then the Betti numbers and the Bass numbers over $\phi$ of any non-zero finitely generated 
$R$-module grow at the same rate, on an exponential scale, as the Betti numbers of $k$ over $R$.
 \end{abstract}


\maketitle

\section{Introduction}

Given an endomorphism $\phi\colon R\to R$ of a commutative Noetherian local ring, each 
$R$-module $M$ defines a module $\up{\phi}M$: it has the same underlying additive group 
as $M$, and $R$ acts on it by the rule $r\cdot m=\phi(r)m$.  We study homological 
properties of $\up{\phi}M$ when $\phi$ is \emph{contracting}; this means that for each $r$ in 
the maximal ideal $\fm$ of $R$ the sequence $(\phi^i(r))_{i\ge1}$ converges to $0$ in the 
$\fm$-adic topology.  

An $R$-module is said to be finite if it is finitely generated over $R$. We prove:

 \begin{theorem}
     \label{ithm:special}
Let $R$ be a local ring and $\phi\colon R\to R$ a contracting endomorphism.

If there exist a finite non-zero $R$-module $M$ and an integer $i\ge1$, such that $\up{\phi^i}M$ has 
finite flat dimension or finite injective dimension, then $R$ is regular.
  \end{theorem}

When the ring $R$ has characteristic $p>0$ and $\phi$ is the Frobenius map, $r\mapsto r^p$, 
the theorem implies that if $\up{\phi}M$ is flat, then $R$ is regular.  We give a second, 
independent argument for this statement. Even when $\up{\phi}M$ is free, it yields a 
substantial strengthening of the classical result of Kunz, \cite[2.1]{Ku}, which treats the 
case $M=R$.  Comparison with other results is given in Remarks \ref{rem:earlier} 
and \ref{rem:evenearlier}.

Other naturally occurring contracting endomorphisms are described in 
Section~\ref{Contractions}.  Here we note that if there is a homomorphism
of rings $R/\fm\to R$, which composed with the natural surjection $R\to R/\fm$ 
gives the identity of $R/\fm$, then the composition of these maps in reverse 
order is a contracting endomorphism of $R$.  Extremal as it is, this 
example captures three motifs that run through the paper:  Contracting
endomorphisms exist in all characteristics; see Example~\ref{ex:existence}.
They exist only for equicharacteristic rings; see Remark~\ref{rem:field}.
When seen through a contracting homomorphism, homological properties of
finite $R$-modules mirror those of $k=R/\fm$.

The proof of Theorem \ref{ithm:special} gives substance to the last point.  It is obtained as a limit case of 
a result that establishes, in precise quantitative terms, that for \emph{every} non-zero 
$M$ over \emph{any} $R$ the (co)homology of $\up{\phi}M$ behaves asymptotically 
as that of $k$.  In particular, we characterize complete intersections in 
parallel with regular rings.

Indeed, let $\ell_R(-)$ denote length over $R$, and define the
\emph{curvature} of $k$ by
\[
\curv_{R}k=\limsup_n \sqrt[n]{\ell_R\Tor^{R}_{n}(k,k)}\,.
\]
It measures, on an exponential scale, the asymptotic rate of growth of 
a minimal free resolution of $k$.  All groups $\Tor^{R}_{n}(k,\up{\phi}M)$ 
and $\Ext_{R}^{n}(k,\up{\phi}M)$ have actions of $R$ induced by the \emph{original} 
action on the additive group shared by $M$ and $\up\phi M$.  The resulting $R$-modules 
are annihilated by the ideal $\phi(\fm)R$, and are finite  when $M$ is.

A version of our main theorem can now be stated as follows: 
 
 \begin{theorem}
     \label{ithm:main}
If $(R,\fm)$ is a local ring, $\phi\colon R\to R$ a contracting endomorphism, and
the ring $R/\phi(\fm)R$ is artinian, then every finite non-zero $R$-module 
$M$ satisfies
\begin{align*}
\limsup_n \sqrt[n]{\ell_{R}\Tor^{R}_{n}(k,\up{\phi}M)} 
=\curv_{R}k=\limsup_n \sqrt[n]{\ell_{R}\Ext_{R}^{n}(k,\up{\phi}M)}\,.
\end{align*}
  \end{theorem}

The preceding results are corollaries of Theorem \ref{thm:endomorphism}, where 
$M$ is a complex with finite homology and the ring $R/\phi(R)$ is not assumed artinian.  
Absent the latter hypothesis, the numbers $\ell_R\Tor^{R}_{n}(k,\up{\phi}M)$ and 
$\ell_R\Ext_{R}^{n}(k,\up{\phi}M)$ need not be finite.

To deal with with this problem we replace lengths with Betti numbers and Bass numbers 
\emph{over the map} $\phi$.  The definition of these numbers, given in Section~\ref{Betti}, 
involves suitable Koszul complexes.  This approach originates in \cite{AIM}, where it was 
developed for bounded complexes with finite homology.  However, that context is too narrow to 
accommodate the proof of Theorem \ref{thm:endomorphism}, even when $M$ is an 
$R$-module.  In Sections \ref{heaven} and \ref{hell} we prove the relevant properties 
of homological invariants over $\phi$, for complexes belonging to appropriate derived 
categories of $R$-modules.  

Section \ref{S:tightcl} can be read independently of the preceding ones.  Using tight closure
methods, 
see \cite{HH}, we give a different proof that a ring $R$ of positive characteristic 
is regular if $\up{\phi^i}M$ is flat for a finite module $M\ne0$ and the Frobenius 
endomorphism~$\phi$. 

%%%%%%%%%%%%%%%%%%%%%%%%%%%%%%%%%%%%%%%%%%%%%%
\section{Asymptotic invariants}
\label{Betti}
%%%%%%%%%%%%%%%%%%%%%%%%%%%%%%%%%%%%%%%%%%%%%%

Let $R$ be commutative ring, $\dcat R$ the derived category of $R$-modules and
$\shift$ the translation functor; $\simeq$ flags isomorphisms in $\dcat R$. Complexes 
carry lower gradings:
 \[
M=\quad\cdots\lra M_{n+1}\xra{\dd^M_{n+1}} M_{n}\xra{\ \dd^M_{n}\ }
M_{n-1}\lra\cdots
 \]
Let $\dcatp R$ (respectively, $\dcatn R$) denote the full subcategory of $\dcat R$ consisting of 
those complexes $M$ for which the $R$-module $\HH_{n}(M)$ is finite for each $n$, and is zero for 
$n\ll 0$ (respectively, $n\gg 0$). Set $\dcatb R =\dcatp R\cap\dcatn R$.  Modules are identified 
with complexes concentrated in degree $0$, and the category of $R$-modules is identified with 
the full subcategory of $\dcat R$ with objects $\{M\in\dcat R\mid \HH_{n}(M)=0\text{ for }n\ne0\}$.

The derived functors of tensor products and of homomorphisms are denoted $-\dtensor R-$ 
and $\Rhom R--$, respectively.  For each integer $n$, we set
 \[
\Tor^R_n(-,-) = \HH_n(-\dtensor R-)
\qquad\text{and}\qquad
\Ext_R^n(-,-) = \HH_{-n}(\Rhom R--)
 \]

  \begin{chunk}
    \label{ch:finite}
Let $\vf\colon R \to S$ be a homomorphism of commutative Noetherian rings.

Complexes of $S$-modules are always viewed as complexes of $R$-modules by restricting 
scalars along $\vf$.  As explained in \cite[1.1]{AIM}, when $M$ and $N$ are complexes of 
$S$-modules the functors $-\dtensor RM$ and $\Rhom R-N$ induce functors
\[
-\dtensor RM \colon \dcat R \to \dcat S\quad\text{and}\quad \Rhom R-N\colon \dcat R \to \dcat S\,,
\]
When $L$ is in $\dcatp R$ with $L\not\simeq0$, $M$ in $\dcatp S$, and $N$ in $\dcatn S$ the following hold:
  \begin{alignat}{4}
    \label{eqn:finite1}
L\dtensor RM&\in\dcatp S
  &&\quad\text{and}\quad
&L\dtensor RM&\not\simeq0
  &\quad\text{when}\quad
 M&\not\simeq0
   \\
    \label{eqn:finite2}
\Rhom RLN&\in\dcatn S
  &&\quad\text{and}\quad
&\Rhom RLN&\not\simeq0
  &\quad\text{when}\quad
N&\not\simeq0
  \end{alignat}
 \end{chunk}  

 \begin{chunk}  
   \label{ch:koszul}
Given a finite subset  $\bsx$ of a commutative ring $S$, let $\koszul\bsx S$ denote the Koszul 
complex on $\bsx$.  
For each complex $M$ of $S$-module, set $\koszul{\bsx}M = \koszul{\bsx}S\otimes_SM$.
The classical isomorphism $\koszul{\bsx}S\cong\shift^{-e}\Hom_S(\koszul{\bsx}S,S)$, where
$e=\operatorname{card}\bsx$, yields an isomorphism
$\koszul{\bsx}M\cong\shift^{-e}\Hom_S(\koszul{\bsx}S,M)$ of complexes of $S$-modules.

Let $(S,\fn,k)$ be a \emph{local ring}; here this means that $S$ is a commutative Noetherian ring with 
unique maximal ideal $\fn$, and $l=S/\fn$ is its residue field. When $\bsx$ is a minimal generating 
set for $\fn$, the complex $\koszul \bsx M$ is independent of the choice of $\bsx$, up to isomorphism, 
so we write $K^{M}$ in place of $\koszul \bsx M$.
   \end{chunk}  

For the rest of the paper, we fix a \emph{local homomorphism} $\vf\colon (R,\fm,k)\to(S,\fn,l)$; that is,
a homomorphism of rings $\vf\colon R\to S$, satisfying $\vf(\fm)\subseteq\fn$.  Set
\[
\edim\vf=\edim(S/\fm S)\,.
 \]

 \begin{chunk}  
   \label{ch:gen}
Let $\bsy$ in $S$ be a \emph{minimal set of generators of $\fn$ modulo $\fm S$}, by which we mean 
that it contains $\edim\vf$ elements and its image in $S/\fm S$ generates the ideal $\fn/\fm S$.

For $M$ in $\dcatp S$ each $S$-module $\Tor^R_n(k,\koszul{\bsy}M)$ is finite, and is 
equal to zero for $n\ll0$; see, \eqref{eqn:finite1}.   It is annihilated by $\fn$, see 
\cite[1.5.6]{AIM}, so it is an $l$-vector space of finite rank.  By definition, the $n$th 
\emph{Betti number of $M$ over $\vf$} is the integer
\[
  \betti n\vf M = \rank_l\Tor^R_n(k,\koszul{\bsy}M)\ge0 \,,
\]
and the \emph{Poincar\'e series} of $M$ over $\vf$ is the formal  Laurent series
  \[
\rP{\vf}{M}=\sum_{n\in\BZ} \betti n\vf M t^n\in\BZ[\![t]\!]\,.
  \]
In case $\vf=\idmap^R$, one gets the usual Betti numbers and Poincar\'e series over $R$. 

When $\bsx$ is a set of generators of $\fn$ containing $q$ elements, the proof of \cite[4.3.1]{AIM} 
(where it is assumed that $M$ is in $\dcatb S$) applies \emph{verbatim} to give an equality
  \begin{equation}
     \label{eqn:poincare}
\rP{\vf}M (1+t)^{q - \edim \vf} =\sum_{n\in\BZ}\rank_l\Tor^R_n(k,\koszul{\bsx}M)t^n\,.
  \end{equation}
Choosing $\bsx$ minimal one sees that $\rP{\vf}M$, and thus $\betti n\vf M$, does not depend on $\bsy$.
   \end{chunk}  

\begin{chunk}
\label{ch:ass}
For $M$ in $\dcatp S$,  the \emph{curvature} and the \emph{complexity} of $M$ over $\vf$ are the numbers
  \begin{align}
    \label{eqn:ass1}
\curv_{\vf}M &= \limsup_n \sqrt[n]{\betti n{\vf}M}
  \\
    \label{eqn:ass2}
\cxy_{\vf}M  &= \inf\left\{d \in\BN \left|
\begin{gathered}
\text{there exists $c\in\BR$ such that}\\
\betti n{\vf}M\leq c n^{d-1}\text{ for all $n\gg0$}
\end{gathered}
\right\}\right.
  \end{align}
In case $\vf=\idmap^R$, we write $\curv_RM$ and $\cxy_{R}M$, respectively.  

When $M$ is in $\dcatb S$ the following inequalities hold, see \cite[7.1.3(5)]{AIM}:
\begin{gather}
    \label{eqn:ass3}
\curv_{\vf}M \leq \curv_{R}k <\infty\quad\text{and}\quad \cxy_{\vf}M \leq \cxy_{R}k\,.
 \end{gather}
If, in addition, the ring $S/\vf(\fm)S$ is artinian, \cite[7.2.3]{AIM} yields
  \begin{equation}
    \label{eqn:ass4}
\curv_\vf M=\limsup_n \sqrt[n]{\ell_{R}\Tor^{R}_{n}(k,M)}\,.
  \end{equation}
  \end{chunk}

\begin{chunk}
\label{ch:ass2}
For $N$ in $\dcatn S$, the $n$th \emph{Bass number $\bass n\vf N$ of $N$ over $\vf$} 
is the integer
\[
  \bass n\vf N = \rank_l\Ext_R^{n-\edim\vf}(k,\koszul{\bsy}N)\ge0\,,
\]
with $\bsy$ as in \ref{ch:gen}, and the \emph{Bass series} of $M$ over $\vf$ is the formal  Laurent series
  \[
\rI{\vf}{M}=\sum_{n\in\BZ} \bass n\vf M t^n\in\BZ[\![t]\!]\,.
  \]

With $\bsx$ as in \ref{ch:gen}, the proof of \cite[4.3.1]{AIM} applies \emph{verbatim} to give an equality
  \begin{equation}
     \label{eqn:bass}
\rI{\vf}M (1+t)^{q - \edim \vf} =\sum_{n\in\BZ}\rank_l\Ext_R^n(k,\koszul{\bsx}M)t^n\,.
  \end{equation}
As above, this implies that $\rI{\vf}M$ and $\bass n\vf N$ are, indeed, invariants of $M$.  

The obvious analogs of \eqref{eqn:ass1} and \eqref{eqn:ass2} define new asymptotic invariants of 
$N$ over $\vf$: its \emph{injective curvature} $\injcurv_{\vf}N$ and its \emph{injective complexity} 
$\injcxy_{\vf}N$.  Furthermore, the analog of \eqref{eqn:ass4} holds, again by \cite[7.2.3]{AIM}.
  \end{chunk}

%%%%%%%%%%%%%%%%%%%%%%%%%%%%%%%%%%%%%%%%%%%%%%
\section{Duality and compositions}
\label{heaven}
%%%%%%%%%%%%%%%%%%%%%%%%%%%%%%%%%%%%%%%%%%%%%%

In this section we study the behavior of complexities and curvatures under  formation of 
Matlis duals and compositions of local homomorphisms.  For expository reasons, we extend
the notation for complexity and curvature.

\begin{chunk}
\label{ch:estimates}
Let $a(t)=\sum_{n=i}^\infty a_nt^n$ be a formal Laurent series with $a_n$ real and non-negative.
 
We set $\curv a(t)=\limsup_n \sqrt[n]{a_{n}}$ and let $\cxy a(t)$ denote the least natural number 
$d$ such that, for some $c\in\BR$ one has $a_{n}\leq cn^{d-1}$ for all $n\gg 0$.  

Let $b(t)=\sum_{n=i}^\infty b_nt^n$ be a Laurent series with $b_n$ real and non-negative.

We write $a(t)\preccurlyeq b(t)$ when $a_n\le b_n$ holds for each $n\in\BZ$; clearly, one then has
\begin{gather}
\label{eqn:cucx1}
 \curv a(t) \leq \curv b(t)  \quad\text{and}\quad \cxy a(t) \leq \cxy b(t)\,.
\end{gather}

The product $a(t)b(t)$ satisfies the following (in)equalities:
\begin{align}
\label{eqn:cucx2}
  \curv(a(t)b(t))&= \max\{\curv a(t), \curv b(t)\} \\
\label{eqn:cucx3}
  \max\{\cxy a(t), \cxy b(t)\} \leq  \cxy(a(t)b(t))& \leq  \cxy a(t) + \cxy b(t) 
\end{align}

Indeed, $\curv a(t)$ is the reciprocal of the radius of convergence of $a(t)$, hence
$\curv(a(t)b(t)) \leq \max\{\curv a(t), \curv b(t)\}$. For the converse, we may assume 
$a_{n}\ne 0$ for some $n$; then $a(t)b(t)\succcurlyeq a_nt^{n}b(t)$ holds, so \eqref{eqn:cucx1} 
yields the inequality below:
\[
\curv(a(t)b(t))\geq \curv(a_nt^{n}b(t))= \curv b(t)\,.
\]
By symmetry, we also have $\curv(a(t)b(t))\geq \curv a(t)$, as desired.

The estimates for $\cxy(a(t)b(t))$ are equally easy to verify.
\end{chunk}

 \begin{proposition}
\label{prop:duality}
If $E$ is an injective hull of\,\ $l$ over $S$ and $M$ a complex in $\dcatn S$, then 
the complex $N=\Hom_{S}(K^{M},E)$ is in $\dcatp S$ and the following equalities hold:
  \begin{gather*}
\rP{\vf}{N}=\rI{\vf}M (1+t)^{\edim S}\,,
  \\
\curv_{\vf}{N}=\injcurv_{\vf}M
  \quad\text{and}\quad
\cxy_{\vf}{N}=\injcxy_{\vf}{M}\,.
  \end{gather*}
\end{proposition}

\begin{proof}
We first show that $\ell_S(\HH_n(K^{M}))$ is finite for each $n$, and is zero for 
$n\gg0$.  

Set $d=\edim S$.  The  filtration $(K^S_{\les p}\otimes_SM)_p$
yields a spectral sequence with
  \[
E_{p,q}^2=\HH_{p}\big(K^{\HH_{q}(M)}\big)
\quad\text{and}\quad {d}_{p,q}^r\colon E_{p,q}^r\to E_{p-r,q+r-1}^r\,.
   \]
The definition of $K^S$ yields $E_{p,q}^2=0$ for $p\le-1$ and for $p\ge(d+1)$.
It follows that $E_{p,q}^r=E_{p,q}^{r+1}$ holds for $r\ge d$, so the spectral sequence 
converges to $\HH_{p+q}(K^{M})$. 

The hypothesis $\HH_{q}(M)=0$ for $q\gg0$ yields $E_{p,q}^2=0$ for $q\gg0$,
which implies $\HH_{n}(K^{M})=0$ for $n\gg0$, due to the convergence of the sequence.  
Moreover, $\HH_p(K^{\HH_{q}(M)})$ is Noetherian along with $\HH_{q}(M)$,
and is annihilated by $\fn$, so each $E_{p,q}^2$ has finite length; the convergence 
of the sequence implies that so does $\HH_n(K^{M})$.

By the injectivity of $E$, for every $n\in\BZ$ there is an isomorphism of $S$-modules
  \[
\HH_n(N)=\HH_n\Hom_{S}(K^M,E)\cong\Hom_{S}(\HH_{-n}(M),E)\,,
  \]
which shows that $\HH_n(N)$ is finite for each $n$ and is zero for $n\ll0$.

Set $e=\edim\vf$, and let $\bsy$ be a minimal generating set of $\fn$ modulo $\fm S$.
From the definitions, \ref{ch:koszul}, and adjunction we get isomorphisms of complexes 
of $S$-modules
  \begin{align*}
\koszul{\bsy}{N} 
  &=\koszul{\bsy}S\otimes_S{\Hom_{S}(K^{M},E)}  \\
  &\cong\Hom_S\big(\shift^{-e}\koszul{\bsy}S,{\Hom_{S}(K^{M},E)}\big)\\
  &\cong\Hom_S\big((\shift^{-e}\koszul{\bsy}S\otimes_SK^{M}),E\big)\\
  &=\Hom_{S}(\shift^{-e}\koszul{\bsy}{K^{M}},E)
\end{align*}
They explain the first one in the following string of isomorphisms in $\dcat S$:
 \begin{align*}
k\dtensor  R\koszul{\bsy}N 
  &\simeq k\dtensor R\big(\shift^{-e} \Hom_{S}(\koszul{\bsy}{K^{M}},E)\big) \\
  &\simeq \Hom_{S}(\shift^{-e}\Rhom Rk{\koszul{\bsy}{K^{M}}},E)\,.
\end{align*}
Th second one holds because $k$ has a resolution by finite free $R$-modules,
while $M$ is in $\dcatn S$ and $E$ is injective.  Since $E$ is an injective envelope of $l$,
we obtain the first and the third isomorphisms of $l$-vector spaces in the string
  \begin{align*}
\Tor_{n}^{R}(k,{\koszul{\bsy}N})
&\cong \Hom_{S}(\Ext^{n-e}_{R}(k,\koszul{\bsy}{K^{M}}),E)\\
&\cong \Hom_{l}(\Ext^{n-e}_{R}(k,\koszul{\bsy}{K^{M}}),\Hom_S(l,E))\\
&\cong \Hom_{l}(\Ext^{n-e}_{R}(k,\koszul{\bsy}{K^{M}}),l)\,.
  \end{align*}
The second isomorphism holds because $\fn$ annihilates $\Ext^{*}_{R}(k,\koszul{\bsy}{K^{M}})$.  
Therefore, $\betti n{\vf}N = \bass n{\vf}{K^{M}}$ holds for each $n$.  From this and
\eqref{eqn:bass}, we get
  \[
\rP{\vf}N = \rI {\vf}{K^{M}} = \rI{\vf}M \cdot (1+t)^{d}\,.
   \]
The formulas for curvature and complexity follow, due to \eqref{eqn:cucx2} and  
\eqref{eqn:cucx3}.
  \end{proof}

For $M=S$, the following result reduces to \cite[9.1.1(1)]{AIM}.

\begin{proposition}
\label{prop:composition}
Let $\rho \colon {R'}\to R$ and $\vf\colon R\to S$ be local homomorphisms. 

For each $L\in\dcatp R$ and $M\in \dcatp S$ there are inequalities:
  \begin{align*}
\curv_{\vf\circ\rho}(L\dtensor RM) &\leq \max\{\curv_{\rho}L,\curv_{\vf}M\}\,,
  \\
\cxy_{\vf\circ\rho}(L\dtensor RM) &\leq \cxy_{\rho}L + \cxy_{\vf}M \,.
  \end{align*}
\end{proposition}

\begin{proof}
Let $\fm'$, $\fm$, and $\fn$ denote the maximal ideals of ${R'}$, $R$, and $S$, respectively. 
Let  ${\bsy'}$ be a minimal generating set of $\fm$ modulo $\fm' R$, let $\bsy$ be one of $\fn$ 
modulo $\fm S$ and set $\bsz = \vf({\bsy'})\sqcup \bsy$.  The isomorphism 
$\fn/\fm S \cong (\fn/\fm' S)/(\fm S/\fm' S)$ implies that $\bsz$ generates $\fn/\fm' S$.
Setting $d=\edim\rho+\edim\vf - \edim(\vf\circ\rho)$, and noticing that $L\dtensor RM$ 
is in $\dcatp T$ by \eqref{eqn:finite1}, we may apply \eqref{eqn:poincare} to obtain
     \begin{equation}
        \label{eqn:dtensor}
\rP{\vf\rho}{L\dtensor RM} (1+t)^{d} =\sum_{n\in\BZ}\rank_l\Tor^R_n(k,\koszul{\bsz}{L\dtensor RM})t^n\,.
  \end{equation}

In the derived category of $S$, the isomorphism
  \[
\koszul\bsz{L\dtensor RM} \simeq\koszul{\bsy'}L\dtensor R\koszul{\bsy} M\,,
  \]
combined with the associativity formula for derived tensor products yields
  \[
\left( {k'}\dtensor {R'}{\koszul{\bsy'} L}\right)\dtensor R {\koszul\bsy M} 
\simeq  {k'}\dtensor {R'}\koszul\bsz{L\dtensor RM}\,.
  \]
This isomorphism gives rise to a standard spectral sequence with
  \[
E_{pq}^2=\Tor_p^R\big(\Tor_q^{R'}({k'},\koszul{\bsy'}L),\koszul{\bsy}M\big)
\implies \Tor_{p+q}^{R'}\left({k'},\koszul{\bsz}{L\dtensor RM}\right).
  \]
The $R$-module $\Tor^{R'}(\koszul{\bsy'}L,{k'})$ is annihilated by $\fm$, so one has 
  \[
\Tor_p^R\big(\Tor_q^{R'}({k'},\koszul{\bsy'}L),\koszul{\bsy}M\big)
  \cong
\Tor_q^{R'}({k'},\koszul{\bsy'}L)\otimes_{k'} \Tor_p^R(k,\koszul{\bsy}M)\,.
  \]
The preceding isomorphism and the convergence of the spectral sequence yield
 \begin{equation}
    \label{eq:ps}
\sum_{n\in\BZ}\rank_l\Tor^R_n(k,\koszul{\bsz}{L\dtensor RM})t^n\preccurlyeq \rP{\rho}L\cdot \rP{\vf}M\,.
 \end{equation}

Combining formulas \eqref{eqn:dtensor} and \eqref{eq:ps}, we get a coefficientwise inequality
 \begin{equation*}
    %\label{eq:ps}
 \rP{\vf\circ\rho}{L\dtensor RM} \cdot (1+t)^d \preccurlyeq \rP{\rho}L\cdot \rP{\vf}M
 \end{equation*}
which, by \eqref{eqn:cucx1}, implies the inequality in the following string:
\begin{align*}
\curv_{\vf\circ\rho}{(L\dtensor RM)}
    &=\curv \rP{\vf\circ\rho}{L\dtensor RM}  \\
    &= \curv\big(\rP{\vf\circ\rho}{L\dtensor RM}\cdot(1+t)^{d}\big) \\
              &\le\curv\big(\rP{\rho}L\cdot \rP{\vf}M\big) \\
              &=\max\{\curv\rP{\rho}L,\curv\rP{\vf}M\}\\
              &=\max\{\curv_{\rho}L,\curv_{\vf}M\}\,.
\end{align*}
The equalities at both ends hold by definition, the other two by \eqref{eqn:cucx2}.

A similar argument, using \eqref{eqn:cucx3}, yields $\cxy_{\vf\circ\rho}(L\dtensor RM)
\leq \cxy_{\rho}L + \cxy_{\vf}M$.
 \end{proof}
 
 %%%%%%%%%%%%%%%%%%%%%%%%%%%%%%%%%%%%%%%%%%%%%%
\section{Homotopical Loewy length}
\label{hell}
%%%%%%%%%%%%%%%%%%%%%%%%%%%%%%%%%%%%%%%%%%%%%%

In this section $(S,\fn,l)$ is a local ring and $M$ a complex of $S$-modules.

We introduce two notions that plays a critical, if behind-the-scenes, role in the proof of our 
main results.  The Loewy length of the complex $M$ is the number
\[
\lol_{S}M = \inf\{i\in\BN\mid \fn^{i}M=0\}\,.
\]
The \emph{homotopical Loewy length} of $M$ is defined in \cite{AIM} to be the number
\[
\dll SM = \inf\{\lol_{S}V\mid M\simeq V\text{ in }\dcat S\}\,.
\]

The proof of the next result is extracted from that of \cite[6.2.2]{AIM}, which provides a more 
precise upper bound for the homotopical Loewy length of $K^{S}$.

\begin{proposition}
\label{prop:dll1}
Every complex $M$ over a local ring $(S,\fn,l)$ satisfies
  \[
\dll S{K^{M}}\le\dll S{K^{S}}<\infty\,.
  \]
   \end{proposition}

\begin{proof}
Set $d=\edim S$.  For all integers $i\gg0$ and all $n\in\BZ$, the subcomplex 
  \[
J^i=\quad 0\to \fn^{i-d}K^S_d\to \fn^{i-d+1}K^S_{d-1}\to\cdots\to\fn^{i-1}K^S_1\to \fn^{i}K^S_{0}\to0
  \]
of $K^S$ satisfies $\HH_n(J^i)=0$, by a well-known result of Serre; see \cite[4.1.6(3)]{barca}.  For
such an $i$, the canonical map $K^S\to K^S/J^i$ is a quasi-isomorphism, so it
represents an isomorphism in $\dcat S$.  Now $\fn^{i}(K^S/J^i)=0$ implies $\dll S{K^{S}}\le i<\infty$.

Set $c=\dll S{K^{S}}$.  Let $\varkappa\colon K^S\xra{\simeq}V$ be an isomorphism in $\dcat S$, with
$\fn^cV=0$, and let $\varepsilon\colon F\xra{\simeq}M$ be a semifree resolution.  The quasi-isomorphisms 
of complexes 
  \begin{gather*}
\xymatrixcolsep{2.5pc}
\xymatrix{
K^M=K^S\otimes_SM \ar@{<-}[r]_-{\simeq}^-{F\otimes_S\varepsilon} & K^S \otimes_SF 
\ar@{->}[r]_-{\simeq}^-{\varkappa\otimes_SF} & V\otimes_SF }
  \end{gather*}
represent an isomorphism $K^M\simeq V\otimes_{S}F$ in $\dcat S$.  It implies $\dll S{K^{M}}\le c$,
since 
  \[
\fn^{c}(V\otimes_{S}F)=(\fn^{c}V)\otimes_{S}F=0\,.
  \qedhere
    \]
  \end{proof}
  
Recall that $V\in\dcat S$ is \emph{formal} if there is an isomorphism $V\simeq \hh V$ in 
$\dcat S$.  

  \begin{remark}
 \label{rem:formal}
If $\hh V$ is projective, then $V$ is formal. 

Indeed, choosing for each $n\in\BZ$ a splitting $\sigma_n\colon\HH_n(V)\to\ZZ_n(V)$
of the canonical surjection $\ZZ_n(V)\to\HH_n(V)$, and composing $\sigma_n$ with 
the inclusion $\ZZ_n(V)\to V_n$, one gets a quasi-isomorphism $\HH(V)\to V$, 
whence an isomorphism $\HH(V)\cong V$ in $\dcat S$.
  \end{remark}

\begin{proposition}
\label{prop:dll2}
Let $(S,\fn,l)$ be a local ring and set $c=\dll S{K^{S}}$.

If $\vf\colon(R,\fm,k)\to(S,\fn,l)$ is a local homomorphism with $\vf(\fm)\subseteq \fn^{c}$,
then for every complex $M$ of $S$-modules the following assertions hold.
  \begin{enumerate}[\quad\rm(1)] %\quad added 10/15/10
  \item
The complex $K^M$ is formal in $\dcat R$.
  \item
For each $L$ in $\dcat R$ there are isomorphisms of graded $l$-vector spaces
\begin{align*}
\Tor_{*}^{R}(L,K^M) &\cong \Tor_{*}^{R}(L,k) \otimes_{k}\HH_{*}(K^M)\,.
\end{align*}
  \item
If $M$ is in $\dcatp S$ and $M\not\simeq 0$, then there are inequalities
\begin{alignat*}{3}
\curv_\vf M&\geq \curv_Rk
 \quad&&\text{and}\quad 
 &\cxy_\vf M&\geq \cxy_Rk\,.
\end{alignat*}
  \item
If $M$ is in $\dcatn S$ and $M\not\simeq 0$, then there are inequalities
\begin{alignat*}{3}
\injcurv_\vf M&\geq \curv_Rk
   \quad&&\text{and}\quad 
&\injcxy_\vf M&\geq \cxy_Rk \,. 
\end{alignat*}
  \item
If $M$ is in $\dcatb S$ and $M\not\simeq 0$, then equalities hold in \rm{(3)} and \emph{\rm(4)}.
  \end{enumerate}
   \end{proposition}

\begin{proof}
(1).  Proposition \ref{prop:dll1} yields in $\dcat S$ an isomorphism $K^M\simeq V$, with 
$\fn^{c}V=0$.  This implies $\fm\cdot V=0$, so $R$ acts on $V$ through 
$k$. Since $k$ is a field, $V$ is formal in $\dcat k$, see Remark \ref{rem:formal}, and hence 
also in $\dcat R$.  

(2).  From (1) we get the first one of the following isomorphisms in $\dcat R$:
\begin{align*}
L \dtensor R{K^M} \simeq L \dtensor R{\hh{K^M}}
        \simeq (L \dtensor R k)\otimes_{k}{\hh{K^M}}\,.
\end{align*}
The second one holds because $\fm\cdot \hh{K^M}=0$. Now pass to homology and use the K\"unneth isomorphism.

(3) and (5). Set $e=\edim S - \edim \vf$ and $h(t)=\sum_{n\in\BZ}\rank_{l}\HH_{n}(K^M)t^n$.

The isomorphism in (2), applied with $L=k$, gives $\rP{\vf}M (1+t)^{e} = \rP Rk \cdot h(t)$.
This explains the middle equality in the following display, where the first and last ones 
hold by definition, while the remaining two come from \eqref{eqn:cucx2}:
\begin{align*}
\curv_\vf M&=\curv \rP{\vf}M  \\
  &= \curv \big(\rP{\vf}M\cdot(1+t)^{e}\big) \\
              &=\curv\big(\rP Rk \cdot h(t)\big) \\
              &=\max\{\curv \rP Rk,\curv h(t)\} \\
              &=\max\{\curv_Rk,\curv h(t)\} \,.
\end{align*}
It remains to note that $\curv h(t)\ge0$ holds, with equality when $M$ is in $\dcatb S$.

A similar argument, using \eqref{eqn:cucx3}, yields the assertions concerning $\cxy_{\vf}M$.

(4) and (5). This follows from (3) and (5), due to Proposition \ref{prop:duality}
  \end{proof}

\section{Contracting endomorphisms}
\label{Contractions}

An endomorphism $\phi\colon R\to R$ of a local ring $(R,\fm,k)$ is said to be
\emph{contracting} if for every $r$ in $\fm$ the sequence $(\phi^i(r))_{i\ge1}$ converges to 
zero in the $\fm$-adic topology of $R$.  Necessary conditions and sufficient conditions for
the existence of such endomorphisms are discussed in Remark \ref{rem:field} and
Example \ref{ex:existence}, respectively.

Now we present the main result of the paper.
 
\begin{theorem}
\label{thm:endomorphism}
Let $\phi\colon R\to R$ a contracting endomorphism of a  local ring $(R,\fm,k)$.

For each $i\geq 1$ and each complex $M$ in $\dcatb R$ with $M\not\simeq 0$ there 
are equalities
 \[
\curv_{\phi^{i}}M =\curv_R k = \injcurv_{\phi^{i}}M\,.
 \]
 \end{theorem}

Some special cases of the theorem are known from earlier work. 

  \begin{remark}
    \label{rem:earlier}
Assume that $M$ is in $\dcatb R$ and $M\not\simeq 0$.  

The equalities in the theorem hold for 
all $i\gg 1$ by \cite[12.1.3]{AIM}.  

When $M$ is a bounded complexes of free $R$-modules, one gets
$\curv_{\phi^{i}}M =\curv_R k$ for all $i\ge1$ by \cite[5.10]{DGI} and \cite[12.1.5]{AIM}.  
When, in addition, the ring $R$ is Gorenstein, \cite[5.11]{DGI} and \cite[12.1.5]{AIM} 
yield $\injcurv_{\phi^{i}}M=\curv_{R}k$ for all $i\ge1$.
  \end{remark}

Extending the notation for modules, we write $\up{\phi}M$ for the complex with the same 
underlying graded abelian group as $M$ and $R$-action given by $r\cdot m=\phi(r)m$.

\begin{proof}[Proof of Theorem \emph{\ref{thm:endomorphism}}]
It suffices to treat the case $i=1$, for $\phi^{i}$ is contracting for each $i\geq 1$.
Moreover, by Proposition~\ref{prop:duality}, it suffices to prove $\curv_{\phi}M=\curv_{R}k$.

Set $M^{(1)}=M$ and for each integer $n\geq 2$ define, inductively, a complex 
\[
M^{(n)} =  {M^{(n-1)}} \dtensor R{\up{\phi}M}
\]
in $\dcat R$, where the action of $R$ on $M^{(n)}$ is obtained by applying \ref{ch:finite}
to $\phi\colon R\to R$, with $L= M^{(n-1)}$ and $M$.  Thus, it is induced by the action on 
the additive group of $\up{\phi}M$, coming from the original action of $R$ on the additive group of~$M$.

We claim that for $n\geq 1$ the following statements hold:
\begin{enumerate}[\quad\rm(1)]
\item[(1${}_n$)] $M^{(n)}$ is in $\dcatp R$ and $\hh{M^{(n)}}\ne 0$.
\item[(2${}_n$)]  $\curv_{\phi^{n}}M^{(n)} \leq \curv_{\phi}M$.
\end{enumerate}

Indeed, both assertions are tautological for $n=1$, so we may assume that they hold for 
some $n\ge1$.  Now \eqref{eqn:finite1} and the induction hypothesis give (1${}_{n+1}$).  
To obtain (2${}_{n+1}$) we use the following relations,  which come from 
Proposition~\ref{prop:composition} applied with $R'=R=S$, $\rho=\phi^{n}$, and $\vf=\phi$, 
and from the induction hypothesis
\begin{align*}
\curv_{\phi^{n+1}}M^{(n+1)} 
        &= \curv_{\phi\circ\phi^{n}}\big( {M^{(n)}} \dtensor R {\up{\phi}M} \big) \\
        &\leq \max\{\curv_{\phi^{n}}M^{(n)},\curv_{\phi}M\}\\
        &\leq \max\{\curv_{\phi}M,\curv_{\phi}M\}\\
        &=\curv_{\phi}M\,.
\end{align*}

Set $c=\dll S{K^{S}}$.  As $\phi$ is contracting, we have $\phi^{s}(\fm)\subseteq\fm^{c}$ for
some integer~$s$.  Applying Proposition~\ref{prop:dll2}(3), assertion (2${}_s$) 
above, and \eqref{eqn:ass3} we now get
\[
\curv_{R}k\leq \curv_{\phi^{s}}M^{(s)} \leq \curv_{\phi}M \leq \curv_{R}k\,.
  \qedhere
\]
  \end{proof}

The notation and hypotheses of the theorem are kept in force in its corollaries.

Part (1) of the first corollary contains Theorem~\ref{ithm:special}, announced 
in the introduction.

\begin{corollary}
\label{cor:special}
For each positive integer $i$ the following hold.
  \begin{enumerate}[\quad\rm(1)] %\quad added 10/15/10
  \item
If $\up{\phi^{i}}M$ is isomorphic in $\dcat R$ to a bounded complex of flat $R$-modules, 
or to a bounded complex of injective $R$-modules, then $R$ is regular.
  \item
If $\curv_{\phi^i}M\le1$ or $\injcurv_{\phi^i}M\le1$ holds, then $R$ is 
complete intersection.
  \end{enumerate}
\end{corollary}

  \begin{remark}
    \label{rem:evenearlier}
Part (1) of the corollary contains Rodicio's generalization of Kunz's Theorem: When $R$ is of 
characteristic $p>0$ and $\phi$ is the Frobenius map, if $\up{\phi^{i}}R$ has finite flat dimension 
for some $i$ then $R$ is regular; see \cite[Thm.~2]{Ro}. 
  \end{remark}

\begin{proof}[Proof of Corollary \emph{\ref{cor:special}}]
(1)  The hypotheses on $\up{\phi^{i}}M$ imply $\betti n{\phi^{i}} M=0$ or $\bass n{\phi^{i}} M=0$
for all $n\gg0$, whence $\curv_{\phi^i}M=0$ or $\injcurv_{\phi^i}M=0$.  The theorem 
then yields $\curv_{R}k=0$, so $R$ is regular by the Auslander-Buchsbaum-Serre Theorem.

(2)  The theorem gives $\curv_{R}k\le1$, so $R$ is complete intersection by \cite[8.2.2]{barca}.
  \end{proof}

The corollary \emph{characterizes} regularity and complete intersection, since it is known 
that the converses of both (1) and (2) hold.  This follows immediately from the precise 
information available on the asymptotic behavior of Betti numbers and Bass numbers over 
contracting endomorphisms of complete intersections.

  \begin{remark}
    \label{rem:special}
When $R$ is complete intersection, \cite[5.3.2]{AIM} yields for each $M$ in $\dcatb R$ 
polynomials $b^M_\pm(t)\in\BQ[t]$ with the same leading term and of degree at most 
$\codim R-1$, such that Betti numbers $\betti n\phi M$ satisfy the equalities
  \[
\betti n\phi M=
\begin{cases}
b^M_+(n)&\text{for all even }n\gg0\,,
  \\
b^M_-(n)&\text{for all odd }n\gg0\,.
\end{cases}
  \]
Furthermore, the Bass numbers $\bass n\phi M$ have a similar property.
  \end{remark}

In \cite{AIM}, a complex $M$ in $\dcatb R$ is said to be \emph{extremal} over 
$\phi$ if it satisfies
 \[
\curv_{\phi}M =\curv_R k
  \quad\text{and}\quad
\cxy_{\phi}M =\cxy_R k  \,.
 \]
The obvious substitutions yield a definition of \emph{injective extremality}.  

Part (b) of the next corollary answers, in the positive, Question \cite[12.2.2]{AIM}.

  \begin{corollary}
    \label{cor:extremal}
Under the following conditions, $M$ is extremal over $\phi^i$ for $i\ge1$:
  \begin{enumerate}[\quad\rm(a)]
    \item
The ring $R$ is not complete intersection.
   \item
The ring $R$ has positive characteristic and $\phi$ is the Frobenius endomorphism.
  \end{enumerate}
  \end{corollary}

  \begin{proof}
Theorem \ref{thm:endomorphism} shows that we need only compare complexities.

Condition (a) implies $\curv_Rk>1$ and $\cxy_Rk=\infty$ by \cite[8.2.2]{barca}, so from
Theorem \ref{thm:endomorphism} we obtain $\curv_{\phi^i}M>1<\injcurv_{\phi^i}M$, 
whence $\cxy_{\phi^i}M =\infty=\cxy_{\phi^i}M$.

Under condition (b), the equalities of complexities are proved in \cite[12.2.4]{AIM}.
  \end{proof}

The restriction in condition (a) is essential:

  \begin{remark}
    \label{rem:extremal}
When $R$ is complete intersection, $M$ is extremal and injectively extremal 
over $\phi^i$ for $i\gg0$, by \cite[12.1.3]{AIM}, but not in general, see \cite[12.1.6]{AIM}.
  \end{remark}

Theorem~\ref{ithm:main} from the introduction is contained in the next corollary.
It follows from Theorem \ref{thm:endomorphism}, formula \eqref{eqn:ass4},
and its analog for injective complexity, see \ref{ch:ass2}.
  
 \begin{corollary}
When $R/\phi(\fm)R$ is artinian the following equalities hold for $i\ge1$:
\begin{xxalignat}{3}
&{\phantom{\square}}
&\limsup_n \sqrt[n]{\ell_{R}\Tor^{R}_{n}(k,\up{\phi^{i}}M)} 
&=\curv_{R}k=\limsup_n \sqrt[n]{\ell_{R}\Ext_{R}^{n}(k,\up{\phi^{i}}M)}\,.
&&\qed
 \end{xxalignat}
   \end{corollary}

In order to apply our results to a given ring $R$, one needs to know that it admits
\emph{some} contracting endomorphism.   Mohan Kumar and Hamid Rahmati
have noticed that such a ring has to be \emph{equicharacteristic}; that is, to satisfy
$\operatorname{char}(k)R=0$.

   \begin{remark}
     \label{rem:field}
If $R$ admits a contracting endomorphism, then it is equicharacteristic.

More precisely, if $\phi\colon R\to R$ is a contracting endomorphism, then the set
   \[
k_0=\{r\in R\mid \phi(r)=r\}
   \]
is a subfield of $R$.  Indeed, it is immediately clear that $k_0$ is a subring of $R$.  
For $r\in k_0\cap\fm$ one has $r\in\bigcap_{j=1}^{\infty}\fm^j=0$.  Thus, each 
non-zero element $r$ of $k_0$ has an inverse in $R$; for every $i\ge1$ it satisfies 
$\phi(r^{-1})=\phi(r)^{-1}=r^{-1}$, so $r^{-1}$ is in $k_0$.
   \end{remark}

Conversely, equicharacteristic rings often have contracting endomorphisms:

   \begin{examples}
     \label{ex:existence}
(1) If $R$ is equicharacteristic and $\operatorname{char}(k)=p>0$, then the 
Frobenius map $r\mapsto r^p$ is a contracting endomorphism.

(2) If $k$ is an arbitrary field, $B$ is a finitely generated subsemigroup of $\BN^n$
for some integer $n$, and $R$ is the localization of $k[B]$ at the maximal ideal
spanned by the positive elements of $B$, then for every integer $q\ge2$ the map
$B\to B$ given by $b\mapsto qb$, induces a contracting endomorphism $R\to R$.

(3) If the canonical map $\varepsilon\colon R\to R/\fm$ admits a left inverse 
homomorphism of rings $\sigma \colon k\to R$, then 
$\sigma\varepsilon\colon R\to R$ is a contracting endomorphism.  

In particular, every equicharacteristic and complete local
ring admits a contracting endomorphism, due to Cohen's Structure Theorem.
   \end{examples}

\section{A proof of a special case for the Frobenius endomorphism} 
\label{S:tightcl}
   
In this section, we give an entirely different proof of a special case of Theorem~\ref{ithm:special}
of the Introduction using tight closure methods in the case where $\phi$ is a 
power of the Frobenius endomorphism of $R$.  Precisely:

\begin{theorem}
\label{thm:Frobenius}
Let $M$ be a finitely generated module over a Noetherian ring $R$ of positive 
prime characteristic $p$ that is supported everywhere on $\Spec(R)$.  Suppose 
that $\phi = F^e$ is an iteration of the Frobenius endomorphism $F$ of $R$ and that 
$\up{\phi}M$ is $R$-flat.  Then $R$ is regular.
\end{theorem}  

We make use of bimodules in the sequel.  If $\phi\colon  R \to S$ is a homomorphism and $W$ is a right $S$-module, then we may define an  $(R,\,S)$-bimodule structure on  $W$ on the abelian group $W$ such that the left $R$-module structure is given by restriction of scalars, that is $\up{\phi}W$, and the right $S$-module structure is the original one.  Thus, $rw =w\phi(r)$ for all $r \in R$ and $w \in W$.  In the sequel, a \emph{$\phi$-bimodule}
means an $(R,\,S)$-bimodule isomorphic to one obtained as above. 

\begin{remark}
\label{rem:tensor}
Let $\phi\colon  R \to S$ and $\psi\colon  S \to T$ be ring homomorphisms, let $M$ be a $\phi$-bimodule and
$N$ a $\psi$-bimodule. The following assertions are clear:

If $M$ is finitely generated as an $S$-module, and $N$ is as a $T$-module, then $M \otimes_S N$ is finitely generated as a $T$-module. 

If $M$ is flat as an $R$-module and $N$ is flat as an {$S$-module}, then $M \otimes_S N$ is flat as an $R$-module.   
\end{remark}

\begin{lemma} 
\label{lem:completion}
Let $\phi\colon  R \to S$ be a local homomorphism, $\wh\phi\colon\wh R\to\wh S$ the induced homomorphism of complete local rings, and $M$ a $\phi$-bimodule.
\begin{enumerate}[\quad\rm(1)]
\item
If $M$ is finitely generated as an $S$-module and flat as an $R$-module, then $M \otimes_S \wh S$ has a structure of $\wh\phi$-bimodule that is finitely generated as an $\wh S$-module and flat as an $\wh R$-module.
\item
If $\Supp_SM=\Spec S$, then $\Supp_{\wh S}(M \otimes_S \wh S)=\Spec\wh S$.
\end{enumerate}
\end{lemma}

\begin{proof}
(1) We have a commutative diagram, with $\iota_R$ and $\iota_S$ the canonical maps:
\begin{equation}
   \label{eq:star}
\begin{gathered}
\xymatrixcolsep{2.5pc}
\xymatrixrowsep{2pc}
\xymatrix{
\wh R \ar@{->}[r]^{\wh {\phi}} & \wh S\\
R \ar@{->}[u]^{\iota_R} \ar@{->}[r]^{\phi} & S\ar@{->}[u]_{\iota_S}                          
}                        
 \end{gathered}
  \end{equation}
Set $\phi'=\iota_S \circ \phi$; this is the same as $\wh{\phi} \circ \iota_R$. Then $M \otimes_S \wh S$ is a right $\wh S$-module, which in turn gets a structure of a $\wh \phi$-bimodule as well as a $\phi'$-bimodule.
Note that $M \otimes_S \wh S$ is finitely generated as an $\wh S$-module and flat as an $R$-module, by Remark~\ref{rem:tensor}, since $M$ is flat over $R$ and $\wh S$ is flat over $S$.

Now, we claim that $M \otimes_S\wh S$ is flat also as an $\wh R$-module. In order to prove this, it suffice to show
that the $\wh R$-module action on $M \otimes_S \wh S$ preserves inclusions of finitely generated $\wh R$-modules. But if  there is a counterexample, then there must be a counterexample involving finite length $\wh R$-modules (cf.\ the Artin-Rees Lemma, the  Krull Intersection Theorem, and Remark~\ref{rem:tensor} showing $M \otimes_S \wh S$ is finitely generated over $\wh S$). And these finite length $\wh R$-modules and the inclusion map in the counterexample must (and trivially) come from the category of $R$-modules via the scalar extension $\iota_R$, which contradicts  the flatness of $M\otimes_S \wh S$ over $R$. (This is  the local flatness criterion.) 

(2) Note that $\Supp_{S}M$ is the set of prime ideals in $S$ containing $\Ann_{S}M$. The desired equality holds because $\Ann_{\wh S}(M \otimes_S\wh S)=(\Ann_{S}M)\wh S$, since $M$ is a finitely generated $S$-module and $\wh S$ is flat over $S$.
\end{proof}

Throughout the remainder of this section,  $R$ will denote a Noetherian ring of prime characteristic $p >0$,  $e$ a positive integer,  and $q$ will denote $p^e$. In this case, an $F^{e}$-bimodule is a (right) $R$-module $M$ with left $R$-module structure given by $rm = mr^q$ for $r \in R$ and $m \in M$. We write $M_{R}$ (respectively, $_{R}M$) for $M$ viewed as a right (respectively, left) $R$-module.

Now Theorem~\ref{thm:Frobenius} may be restated as follows:

\begin{theorem}
\label{thm:FrobeniusAlt}
Let $R$ be a Noetherian commutative ring of prime characteristic $p$. Assume there exists an $F^e$-bimodule $M$ with $e \ge 1$ such that $\Supp(M) = \Spec(R)$, the right $R$-module $M_R$ is finitely generated $R$-module and the left $R$-module $_R M$ is flat.  Then $R$ is regular.  
\end{theorem}

\begin{proof}
We first note that the result reduces at once to the local case. Henceforth, we may assume without loss of generality that $(R,\fm)$ is local.  The proof proceeds in three steps.  We first show that if $R$ is Cohen-Macaulay, then $R$ is regular. Next we show that $R$ must be a domain. Finally, we use tight closure theory to prove that $R$ is, in fact Cohen-Macaulay and hence regular; 
and we achieve this by reducing to the case where $R$ is complete.

As usual, if $I$ is an ideal of $R$,  $I^{[q]}$  denotes the ideal $(f^q: f \in I)R$,  which is the expansion
of $I$ under $F^e\colon R \to R$.   

For every positive integer $n$, denote $M^{(n)}: = M^{\otimes n}$. That is, form  $M^{(n+1)}$ recursively as  $M \otimes_R M^{(n)}$. Note that for all $n$, $M^{(n)}$ is naturally an $F^{ne}$-bimodule that is finitely generated as a right $R$-module and flat as a left $R$-module, by Remark~\ref{rem:tensor}. In fact, $M^{(n)}$ is automatically faithfully flat as a left $R$-module, since $(R/\fm) \otimes_R M^{(n)} \cong M^{(n)}/(M^{(n)}\fm^{[q^n]})$ as right $R$-modules, and the latter module is nonzero by Nakayama's lemma.  

First, we show that if $(R,\fm)$ is Cohen-Macaulay, then $R$ is regular.  In this case  $M$ is a (possibly big)
Cohen-Macaulay left $R$-module; hence $M$ is a (small) Cohn-Macaulay right $R$-module. The same is true for $M^{(n)}$. Let $\bsx$ be a full system of parameters for $R$ (hence, an $R$-regular sequence). There exists $n$
big enough such that $\fm M^{(n)} \subseteq M^{(n)}(\bsx)$. (Note that $\fm^{[q^{n}]}
\subseteq (\bsx)$ for $n \gg 0$.) Replace $M$ by $M^{(n)}$  so that we may assume
$\fm M \subseteq M(\bsx)$ in the remainder of this part of the argument.

Then by considering the complex $M \otimes K[\bsx; R]$, we see $M/(M(\bsx))$ has finite flat dimension as a left $R$-module. Moreover, in light of $\fm M \subseteq M(\bsx)$, we see that $M/(M(\bsx))$ is a non-zero vector space over $R/\fm$, through its structure as a left $R$-module. Thus $R/\fm$ has finite flat dimension hence $R$ is regular, as claimed. 

Next, we show that $(R,\fm)$ must be a domain: Let $a \in R \setminus \{0\}$ and consider the exact sequence 
\[
\xymatrixcolsep{1pc}
\xymatrix{
0 \ar@{->}[r] &  I \ar@{->}[r] & R \ar@{->}[r]^{a} & R},
\]
where $I := (0:_R a)$. Apply $\otimes_R M$ and use the left flatness of $_RM$ to get
\[
\{x \in M \,|\, ax = 0\} = IM = MI^{[q]}.
\]
However, we see directly that
\[
MI\subseteq \{x \in M \,|\, ax = 0\}.
\]
Thus $MI^{[q]} = MI$, and so $MI^2 = MI$. By Nakayama's lemma, we see that 
$MI = 0$. Hence $IM = 0$ or $I \otimes_R M = 0$. It follows that $I=0$
by the faithful flatness of $M$ as a left $R$-module. Therefore, $R$
is automatically a domain as claimed. 

Henceforth, we assume $(R,\fm)$ is local (hence a domain), and prove
the theorem by induction on $\dim(R)$. Being a domain, when
$\dim(R)=0$, $R$ is a field and hence regular. Thus we may assume $\dim(R) \ge 1$.   

Let $F^e_R\colon  R \to R$ be the Frobenius endomorphism, which is local. Noting that $\wh{F^e_R}$ is the Frobenius endomorphism $F^e_{\wh R}$ of $\wh R$, and that $R$ is regular when $\wh R$ is regular, one may assume, by Lemma~\ref{lem:completion}, that $R$ is complete local. 

Therefore, we further assume $(R,\fm)$ is a complete domain. By the result proved at the outset, it suffices to show $R$ is Cohen-Macaulay. In fact, we are going to show $R$ is \emph{weakly F-regular} (that is, $I^* = I$ for all ideals of $R$). To this end, it is enough to show $I = I^*$ for all $\fm$-primary ideals $I$ of $R$. By way of contradiction, suppose there exists $x \in R\setminus I$ such that $x \in I^*$ for some $\fm$-primary ideal $I$. Set $J:= (I, x)$. By choosing $x$ to be a socle element, we may assume $\ell(J/I) = 1$.  

Observe that, for every $P \in\Spec(R) \setminus \{\fm\}$, the induction hypothesis applies to $R_P$ and $M_P$, which shows $R_P$ is regular. In other words, $R$ is an isolated singularity. Thus, as $R$ is a complete (hence excellent) domain, the \emph{test ideal} of $R$, denoted $\tau$, is $\fm$-primary:  see \cite[Theorem 6.20]{HH}.

In what follows, to indicate that we are taking the number of generators, or length, or the annihilator, of a given module viewing it as a right module, we include the superscript $\,^r\,$ in the notation. For instance, $\nu^r(M)$ is the minimal numbers of generators of $M_R$. Clearly, we have  
\[
\nu^r(M) = \ell^r(M/M\fm).
\]
As $R$ is a domain, $\dim(R) \ge 1$ and $M$ is torsion-free (as a left and hence a right $R$-module), we have $M\fm \neq 0$. Thus, Nakayama's implies $M\fm \supsetneq M\fm^2 \supseteq M\fm^{[q]}$. Next, setting
\begin{align*}
f &= \ell^r(k \otimes_R M)= \ell^r(M/\fm M) =
\ell^r(M/M\fm^{[q]})
\quad \text{and} \\
g & = \nu^r(M) = \ell^r(M/M\fm),
\end{align*}
the  argument in the above paragraph implies  
\begin{equation}
  \label{eq:dagger}
f = \ell^r(M/\fm M) = \ell^r(M/M\fm^{[q]}) > \ell^r(M/M\fm) = \nu^r(M) = g. 
  \end{equation}
(Note the strict inequality in \eqref{eq:dagger}.) 

We are going to study $\ell^r(JM^{(n)}/IM^{(n)}) =
\ell^r(M^{(n)}J^{[q^n]}/M^{(n)}I^{[q^n]})$ and get a contradiction.  
As $_R M$ and hence $M^{(n)}$ are flat as left $R$-modules, for all
$n$ we get 
\begin{equation}
  \label{eq:hash}
  \ell^r(JM^{(n)}/IM^{(n)}) = \ell^r(k \otimes_R M^{(n)}) = f^n \,.
  \end{equation}

Also notice that $\nu^r(M^{(n)}) \le g^n$ for all $n$. (In fact,
$\nu^r(M^{(n)}) = g^n$ for all $n$.)

Now let us study $\ell^r(JM^{(n)}/IM^{(n)}) =
\ell^r(M^{(n)}J^{[q^n]}/M^{(n)}I^{[q^n]})$ via the uniform annihilating
property of $\tau$. As $J^{[q^n]} = (I^{[q^n]}, x^{q^n})$, we see
that the minimal numbers of generator of $M^{(n)}J^{[q^n]}/M^{(n)}I^{[q^n]}$
as a right $R$-module satisfies
\[
\nu^r(M^{(n)}J^{[q^n]}/M^{(n)}I^{[q^n]}) \le \nu^r(M^{(n)}) \le g^n
\]
for all $n$. Moreover, we have 
\[
J^{[q^n]} \tau \supseteq I^{[q^n]}, 
\]
which implies that $M^{(n)}J^{[q^n]}/M^{(n)}I^{[q^n]}$, as a
right $R$-module, is killed by $\tau$ (which has been observed to be
$\fm$-primary) for all $n$. Therefore, we see
\[
\ell^r(M^{(n)}J^{[q^n]}/M^{(n)}I^{[q^n]}) \le
\nu^r(M^{(n)}J^{[q^n]}/M^{(n)}I^{[q^n]}) \ell(R/\tau) \le g^n \ell(R/\tau)
\]
for all $n$. Consequently, we get
\begin{equation}
  \label{eq:hash2}
\ell^r(JM^{(n)}/IM^{(n)}) = \ell^r(M^{(n)}J^{[q^n]}/M^{(n)}I^{[q^n]}) \le g^n \ell(R/\tau)
  \end{equation}
for all $n$. Finally, as $g < f$ (see \eqref{eq:dagger}), we must have
\[
\ell^r(JM^{(n)}/IM^{(n)}) \le g^n \ell(R/\tau) < f^n = \ell^r(JM^{(n)}/IM^{(n)})
\]
for all sufficiently large $n$, which is a contradiction. (In other
words, \eqref{eq:hash} and \eqref{eq:hash2}  contradict each other.)

Thus $R$ is Cohen-Macaulay and therefore $R$ is regular.
 \end{proof}

\end{document}